\documentclass[12pt]{article}
\usepackage{amsmath}
\usepackage{latexsym}
\usepackage{amssymb}
\newtheorem{thm}{Theorem}[section]
\newtheorem{la}[thm]{Lemma}
\newtheorem{Defn}[thm]{Definition}
\newtheorem{Remark}[thm]{Remark}
\newtheorem{Note}[thm]{Note}

\newtheorem{Example}[thm]{Example}
\newtheorem{Examples}[thm]{Examples}
\newtheorem{Problems}[thm]{Problems}

\newtheorem{Problem}[thm]{Problem}
\newtheorem{Convention}[thm]{Convention}
\newtheorem{Number}[thm]{\!\!}
\newenvironment{defn}{\begin{Defn}\rm}{\end{Defn}}

\newenvironment{example}{\begin{Example}\rm}{\end{Example}}

\newenvironment{rem}{\begin{Remark}\rm}{\end{Remark}}

\newenvironment{proof}{{\noindent\bf Proof.}}%
                  {\nopagebreak\hspace*{\fill}$\Box$\medskip\medskip\par}   
\newcommand{\Punkt}{\nopagebreak\hspace*{\fill}$\Box$}
\newcommand{\wb}{\overline}

\newcommand{\mto}{\mapsto}

\newcommand{\N}{{\mathbb N}}
\newcommand{\R}{{\mathbb R}}
\newcommand{\bO}{{\mathbb O}}

\newcommand{\K}{{\mathbb K}}

\newcommand{\sub}{\subseteq}

\DeclareMathOperator{\im}{im}

\DeclareMathOperator{\Supp}{supp}

\begin{document}
\renewcommand{\thefootnote}{\fnsymbol{footnote}}
\begin{center}
{\Large\bf
Ultrametric and Non-Locally Convex\\[.6mm]
Analogues
of the General Curve Lemma\\[2.9mm]
of Convenient Differential Calculus}\\[6mm]
{\bf Helge Gl\"{o}ckner\,\footnote{These investigations were
supported by the German Research Foundation
(DFG), project 436 RUS 17/67/05.}}\vspace{2.5mm}
\end{center}
\renewcommand{\thefootnote}{\arabic{footnote}}
\setcounter{footnote}{0}
\begin{abstract}\vspace{1mm}
\hspace*{-7.2 mm}
The General Curve Lemma is a
tool of Infinite-Dimensional Analysis,
which enables refined studies of
differentiability properties of maps
between real locally convex spaces.
In this article, we generalize the General
Curve Lemma in two ways:
First, we remove the condition of
local convexity in the real case.
Second, we adapt the lemma to the case
of curves in topological
vector spaces over ultrametric
fields.\vspace{2.6mm}
\end{abstract}
{\footnotesize {\em Classification}:
26E15, 
26E20, 
26E30, 
45T20, 
46A16, 
46S10\\[1.7mm] 
}\vspace{4mm}
\begin{center}
{\Large\bf Introduction}\vspace{.3mm}
\end{center}
The General Curve Lemma
(as in \cite[Proposition~4.2.15]{FaK} or
\cite[Lemma~12.2]{KaM})
is a powerful
tool for the study of finite order differentiability
properties of mappings between real locally
convex spaces in the Convenient Differential Calculus
of Fr\"{o}licher, Kriegl and Michor
(see \cite[\S\,4.3]{FaK} and \cite[\S\,12]{KaM}).
It allows pieces of a (suitable) given sequence of smooth
curves to be combined to a single smooth curve,
which runs through all of the pieces in finite time.
The goal of this paper is to
extend the General Curve Lemma
to curves in not necessarily locally convex
real topological vector spaces,
and to curves in topological vector spaces
over an ultrametric field.\\[3mm]
Our studies are based on the differential
calculus of smooth and $C^k$-maps
between open subsets of topological vector
spaces over a topological
field developed in~\cite{BGN},
which has by now been applied
to a variety of questions in
Differential Geometry~\cite{Ber},
Lie Theory (\cite{NOA}, \cite{ANA},
\cite{ZOO}) and Dynamical Systems
(see \cite{SUR} for a survey).
We recall that this approach generalizes
traditional concepts:
In particular, a map between open subsets
of real locally convex spaces
is $C^k$ in the sense of~\cite{BGN}
if and only if it is a
Keller $C^k_c$-map (see \cite{BGN}).
Furthermore, it is known (see \cite[Theorem~2.1]{COM})
that a map between open subsets of finite-dimensional
vector spaces over a complete ultrametric field
is $C^k$ in the sense of~\cite{BGN}
if and only if it is
a $C^k$-map in the usual sense of Non-Archimedian
Analysis (as in~\cite[\S\,84]{Sch} and \cite{DSm}).
The definition of $C^1$-maps in~\cite{BGN}
is also similar in spirit
to an earlier definition used in \cite{Ld1}
and~\cite{Ld2}.\\[3mm]
Our General Curve Lemma in the real case
(Theorem~\ref{thmrealcu})
closely resembles its classical counterpart
for curves in real locally convex spaces.
It subsumes:\\[3mm]
\noindent
{\bf Real Case of General Curve Lemma.}
\emph{Let $E$ be a real topo\-lo\-gical
vector space and $(s_n)_{n\in \N}$
as well as $(r_n)_{n\in \N}$
be sequences of positive reals such that
$\sum_{n=1}^\infty s_n<\infty$
and $r_n\geq s_n+\frac{2}{n^2}$
for each $n\in \N$.
Let $(\gamma_n)_{n\in \N}$ be a sequence of smooth
maps $\gamma_n\colon [-r_n,r_n] \to E$
which become small sufficiently
fast $($in the sense
made precise in Theorem~{\rm\ref{thmrealcu})}.
Then there exists a smooth curve $\gamma\colon \R\to E$
and a convergent sequence $(t_n)_{n\in \N}$
of real numbers such that
$\gamma(t_n+t)=\gamma_n(t)$
for all $n\in \N$ and $t\in [-s_n,s_n]$.}\\[3mm]
If $(\K,|.|)$ is an ultrametric field,
we obtain a variant of the General
Curve Lemma (Theorem~\ref{thmgencu})
which subsumes the following result:\\[3mm]
{\bf Ultrametric General Curve Lemma.}
\emph{Let $E$ be a topological vector space
over an ultrametric field~$(\K,|.|)$,
and $\bO:=\{x\in \K\colon |x|\leq 1\}$.
Let $\rho\in \K^\times$ with $|\rho|<1$
and $(\gamma_n)_{n\in \N}$ be a sequence
of maps $\gamma_n\in BC^\infty(\rho^n\bO, E)$
which become small sufficiently fast
$($in the sense made precise
in Theorem~{\rm \ref{thmgencu})}.
Then there exists a smooth map
$\gamma\colon \K\to E$
such that
$\gamma(\rho^{n-1}+t)=
\gamma_n(t)$\linebreak
for all $n\in \N$ and $t\in \rho^n\bO$.}\\[3mm]
The preceding results are useful for the study of
$k$~times H\"{o}lder differentiable maps of H\"{o}lder
exponent $\sigma \in \;]0,1]$ ($C^{k,\sigma}$-maps,
for short),
as introduced in
\cite{COM} and (for $\sigma=1$) in \cite{IM2}.
As shown in \cite{COM},
our General Curve Lemmas imply a characterization
of $C^{k,\sigma}$-maps on
metrizable spaces:\\[3mm]
{\bf Theorem.}
\emph{Let $\K$ be
$\R$ or an ultrametric field.
Let $E$ and $F$ be topological $\K$-vector spaces
and $f\colon U\to F$ be a map,
defined on an open subset $U\sub E$.
Let $k\in \N_0$ and $\sigma\in \;]0,1]$.
If $E$ is metrizable,
then $f$ is $C^{\,k,\sigma}$
if and only if $f\circ \gamma\colon \K^{k+1}\to F$
is $C^{\,k,\sigma}$, for each
$C^\infty$-map
$\gamma \colon \K^{k+1}\to U$.}\\[3mm]
It would be nice to know
whether smooth maps on $\K^{k+1}$
can be replaced by smooth maps
of a single variable here,
as in Boman's classical results
concerning the real finite-dimensional case~\cite{Bom}
and their infinite-dimensional generalizations~\cite{KaM}.
The author undertook some steps in this direction
jointly with S.\,V. Ludkovsky
(cf.\ also Ludkovsky's preprint~\cite{Lud}).
Our versions of the General Curve Lemma
were created in connection with this question.\\[3mm]
We mention that an analogue of the preceding theorem
for $C^k$-maps can\linebreak
already be found in \cite[Theorem~12.4]{BGN},
where it was proved with the help of variants
of the Special Curve Lemma (Lemma~11.1 and
Lemma~11.2 in~\cite{BGN}).\\[3mm]
Our versions of the General
Curve Lemma are more difficult to prove than
the classical lemma (as reflected by the length of this
text), because it does not suffice to prove
merely the existence and continuity of derivatives (of all orders)
for~$\gamma$.
Instead, to establish smoothness
of~$\gamma$, one has to prove existence of
continuous extensions to higher
difference quotient maps, which is a much
more cumbersome task. To keep the effort manageable,
our strategy is to manufacture,
in a first step, certain
smooth curves $\eta_n\colon \K\to E$
with pairwise disjoint supports
from the given curves~$\gamma_n$.
In a second step, we then show that
$\gamma:=\sum_{n=1}^\infty\eta_n$
converges in $BC^\infty(\K,E)$.
To prove convergence of this series,
we introduce a notion of ``absolute
convergence'' for series in
general topological vector spaces
(Definition~\ref{defabscon}),
the topology of which need not arise
from a family of continuous seminorms.
In~\cite{IM2}, so-called ``gauges''
have already been used as
a substitute for
continuous seminorms
(cf.\ \cite{Jar} for the real case).
To define absolute convergence
of series in general topological vector spaces,
we introduce ``calibrations''
as a further generalization of continuous
seminorms (Definition~\ref{defcalib}).
These are sequences of gauges which are
pairwise related by a certain substitute for
the triangle inequality.
\section{Preliminaries, notation and basic facts}\label{secperlim}
%
%
%
%
In this section, we set up terminology
and notation.
We also compile various
basic facts, for later use.
These are easy to take on faith,
and we recommend to skip
the proofs (given in Appendix~A),
which are not difficult.\\[3mm]
All topological fields occurring in this article
are assumed Hausdorff and non-discrete.
A field~$\K$, equipped with
an absolute value $|.|\colon \K\to [0,\infty[$
defining a non-discrete topology on~$\K$
is called a \emph{valued field}.
An \emph{ultrametric field}
is a valued field $(\K,|.|)$
whose absolute value satisfies
the ultrametric inequality,
$|x+y|\leq\max\{|x|,|y|\}$
for all $x,y\in \K$. If $(E,\|.\|)$
is a normed space over a valued
field, $r>0$ and $x\in E$,
we define
$B_r^E(x):=\{y\in E\colon
\|y-x\|<r\}$
and
$\wb{B}_r^E(x):=\{y\in E\colon
\|y-x\|\leq r\}$.
Recall that if $\K$ is an ultrametric field,
then $\wb{B}_r^{\,\K}(x)$ and $B_r^{\,\K}(x)$
are both open and closed
(this will useful for piecewise
definitions of maps).
Furthermore,
the ultrametric inequality implies that
%
\begin{equation}\label{winner}
|x+y|\;=\;|x|\quad\mbox{for all $x,y\in \K$
such that $|y|<|x|$.}
\end{equation}
All topological vector spaces over
topological fields are assumed Hausdorff.
As usual, $\N:=\{1,2,\ldots\}$ and $\N_0:=\N\cup\{0\}$.\\[2.5mm]
A differential calculus of $C^k$-maps between
subsets of ultrametric fields was developed
in~\cite{Sch}. It makes sense just as well
for maps into topological vector spaces
over general topological fields
(cf.\ \cite[\S\,6]{BGN} for open domains),
and will be used in this form here.
The approach can be generalized
to a differential calculus of $C^k$-maps
between open subsets of topological
vector spaces~\cite{BGN}.
Compare \cite{Ld1}, \cite{Ld2}
for an earlier approach to infinite-dimensional
calculus over ultrametric fields,\footnote{See
also the maps called $C^n$ (in contrast to $C^{[n]}$)
in~\cite{Lud}.}
which however is not equivalent to
ours, at least when
applied to local fields
of positive characteristic~\cite{COM}.
We only give the definition
of $C^k$-maps on subsets of~$\K$ here,
following the notational conventions
from \cite{BGN} (rather~than~\cite{Sch}).
%
%
%
\begin{defn}\label{def1}
Let $\K$ be a topological field,
$U\sub \K$ be a non-empty subset without
isolated points,
and $\gamma\colon U\to E$ be a map
to a topological $\K$-vector space~$E$.
The map $\gamma$ is said to be $C^0_\K$
if it is continuous; in this case, we
set $\gamma^{<0>}:=\gamma$.
We call $\gamma$ a \emph{$C^1_\K$-map}
if it is continuous and if there exists
a continuous map $\gamma^{<1>}\colon U\times U\to E$
such that
\[
\gamma^{<1>}(x_0,x_1)\; =\; \frac{\gamma(x_1)-\gamma(x_0)}{x_1-x_0}\quad
\mbox{for all $x_0,x_1\in U$ such that $x_0\not=x_1$.}
\]
Recursively, having defined
$C^j_\K$-maps and associated maps
$\gamma^{<j>}\colon U^{j+1}\to E$
for $j=0,\ldots, k-1$ for some $k\in \N$,
we call $\gamma$ a \emph{$C^k_\K$-map}
if it is $C^{k-1}_\K$
and there is a continuous map
$\gamma^{<k>}\colon U^{k+1}\to E$
such that
\[
\gamma^{<k>}(x_0,x_1,\ldots, x_k)
\;=\;
\frac{\gamma^{<k-1>}(x_k,x_1,\ldots, x_{k-1})
-\gamma^{<k-1>}(x_0,x_1,\ldots, x_{k-1})}{x_k-x_0}
\]
for all $(x_0,\ldots,x_k) \in U^{k+1}$
such that $x_0\not=x_k$.
The map $\gamma$ is $C^\infty_\K$
(or \emph{smooth})
if it is $C^k_\K$
for each $k\in \N_0$.
If $\K$ is understood, we write
$C^k$ instead of $C^k_\K$.
We let $C^k(U,E)$ be the set
of all $C^k$-maps $U\to E$.
Then $C^k(U,E)$
is a vector subspace of $E^U$.
\end{defn}
Here $\gamma^{<k>}$ is uniquely determined,
and $\gamma^{<k>}$ is symmetric in its $k+1$
variables.
Also $k!\, \gamma^{<k>}(x,\ldots,x)
=\frac{d^k\gamma}{dx^k}(x)=:\gamma^{(k)}(x)$,
for all $x\in U$ (cf.\ \cite[\S\,29]{Sch}
and \cite[Proposition~6.2]{BGN}).
Let $U^{>k<}$ be the set of all $(x_0,\ldots, x_k)\in
U^{k+1}$ such that $x_i\not=x_j$
for all $i\not=j$. Then $U^{>k<}$
is dense in $U^{k+1}$, which will be useful later.
\begin{defn}
Let $E$ be a topological
vector space over a topological field~$\K$.
\begin{itemize}
\item[(a)]
A subset $A\sub E$ is called
\emph{bounded} if, for each $0$-neighbourhood
$U\sub E$, there exists a $0$-neighbourhood
$V\sub\K$ such that $V A\sub U$.
\item[(b)]
If~$X$ is a topological space,
then $BC(X,E)$ denotes the set
of all continuous maps $\gamma\colon X\to E$
whose image $\gamma(X)$ is bounded in~$E$.
Clearly $BC(X,E)$ is a vector subspace
of $E^X$. We equip $BC(X,E)$
with the topology of uniform convergence.
\item[(c)]
If $k\in \N_0\cup\{\infty\}$ and
$U\sub \K$ is a non-empty subset
without isolated points,
we let
$BC^k(U,E)$ be the space of all $C^k$-maps
$\gamma\colon U\to E$ such that
$\gamma^{<j>}\in BC(U^{j+1},E)$
for all $j\in \N_0$ such that $j\leq k$.
We equip $BC^k(U,E)$ with the initial
topology with respect to the sequence
of mappings
$BC^k(U,E)\to BC(U^{j+1},E)$,
$\gamma\mto \gamma^{<j>}$
(for $j\in \N_0$, $j\leq k$).
\end{itemize}
\end{defn}
Recall that a topological
vector space over a topological
field~$\K$ is called
\emph{complete}
if each Cauchy net
converges.
%
%
\begin{la}\label{completeBC}
Let $\K$ be a topological field,
$X$ be a topological space,
$U\sub \K$ be a non-empty subset
without isolated points, and $E$
be a topological $\K$-vector space.
Then the following holds:
\begin{itemize}
\item[\rm(a)]
$BC(X,E)$
is a topological $\K$-vector space.
\item[\rm (b)]
If $E$ is complete, then also
$BC(X,E)$ is complete.
\item[\rm (c)]
For each $k\in \N_0\cup\{\infty\}$,
the map $\theta\colon BC^k(U,E)\to\prod_j BC(U^{j+1},E)$,
$\gamma\mto (\gamma^{<j>})_j$ $($where $j\in \N_0$
such that $j\leq k)$ is linear,
a topological embedding and has closed
image.
\item[\rm (d)]
$BC^k(U,E)$
is a topological $\K$-vector space,
for each $k\in \N_0\cup\{\infty\}$.
If~$E$ is complete, then also
$BC^k(U,E)$ is complete.
\end{itemize}
\end{la}
A topological vector space over a valued
field is called \emph{polynormed}
if its vector topology can be defined
by a family of seminorms.
As a replacement
for seminorms
when dealing with
non-polynormed topological vector spaces
over a valued field,
the more general concept of a \emph{gauge}
was introduced in~\cite{IM2}
(cf.\ \cite[\S\,6.3]{Jar} for the real case).
Using gauges, it is easy to define
Lipschitz continuous,
Lipschitz differentiable,
strictly differentiable,
totally differentiable
and similar maps
between arbitrary topological
$\K$-vector spaces (\cite{IM2}, \cite{COM}).
We shall slightly generalize the concept of
a gauge from~\cite{IM2} here,
because this will simplify the presentation
(see Remarks~\ref{oldnew1} and~\ref{oldnew2}).
\begin{defn}
Let $E$ be a topological
vector space over a valued field
$(\K,|.|)$. A \emph{gauge on~$E$}
is a
map $q \colon E\to [0,\infty[$
(also written $\|.\|_q :=q$)
satisfying $q(t x) = |t|q(x)$ for all $t\in \K$
and $x\in E$,
and such that $B_r^q(0):=q^{-1}([0,r[)$
is a $0$-neighbourhood in~$E$, for each $r>0$.
\end{defn}
Note that each gauge is continuous at~$0$.
Sums of gauges and non-negative multiples $rq$ of gauges
are gauges.
%
\begin{rem}\label{oldnew1}
In \cite{IM2},
only upper semicontinuous gauges $q\colon E\to [0,\infty[$
were considered.
Thus, the stronger requirement was made
that $B_r^q(0)$ is open in~$E$,
for each $r>0$. By the next
two remarks, it does not matter
for many purposes
whether the weaker or the stronger
definition is used.
\end{rem}
%
%
\begin{rem}\label{minkow}
Typical examples
of
gauges are Minkowski functionals~$\mu_U$
of balanced, open $0$-neighbourhoods~$U$
in a topological vector space~$E$
over a valued field~$\K$
(see \cite[Remark\,1.21]{IM2}).
These are upper semicontinuous.\linebreak
Here $U\sub E$ is called
\emph{balanced} if $tU\sub U$
for all $t\in \K$ such that $|t|\leq 1$.
The Minkowski functional
is
$\mu_U\colon E\to [0,\infty[$,
$x\mto \inf\{|t|\colon \mbox{$t\in \K^\times$
with $x\in tU$}\}$.
\end{rem}
\begin{rem}
If $q$ is a gauge on~$E$,
then $q\leq \mu_U$
for the Minkowski functional
of some balanced, open
$0$-neighbourhood~$U$.
In fact, we can take
any balanced, open
$0$-neighbourhood $U\sub E$ such that
$U\sub B_1^q(0)$. Given $x\in E$ and $t_n\in \K^\times$
such that $|t_n|\to\mu_U(x)$,
we then have $x\in t_nU\sub B^q_{|t_n|}(0)$
for each~$n$
and thus $q(x)<|t_n|$,
from which $q(x)\leq \mu_U(x)$
follows by letting $n\to\infty$.
\end{rem}
\begin{example}
Given $r\in \;]0,1]$,
a gauge $q\colon E\to [0,\infty[$
is called an $r$-seminorm if
$q(x+y)^r\leq q(x)^r+q(y)^r$
for all $x,y\in E$.
If, furthermore,
$q(x)=0$ if and only if $x=0$,
then~$q$ is called an $r$-norm
(cf.\ \cite[\S\,6.3]{Jar}
for the real case).
For examples of $r$-normed spaces over~$\R$
and more general non-locally convex
real topological vector spaces,
the reader is referred to \cite[\S\,6.10]{Jar}
and~\cite{Kal}.
For $\K$ a valued field,
the simplest
examples are the spaces
$\ell^p(\K)$
of all $x=(x_n)_{n\in\N}
\in \K^\N$
such that $\|x\|_p:=\sqrt[p]{\sum_{n=1}^\infty|x_n|^p}<\infty$,
for $p\in \;]0,1]$.
Then $\|.\|_p$
is a $p$-norm
on $\ell^p(\K)$ defining a Hausdorff
vector topology on this space.
\end{example}
Note that the triangle
inequality need not hold for gauges.
The following lemma (see \cite[Lemma~1.29]{IM2})
provides a certain substitute.
%
%
\begin{la}\label{substitut}
If $E$ is a topological vector
space over a valued field~$\K$
and $U,V\sub E$ are balanced open $0$-neighborhoods
such that $V+V\sub U$, then
\[
\mu_U(x+y)\;\leq\; \max\{\mu_V(x),\mu_V(y)\}\quad
\mbox{for all $\,x,y\in E$.}
\]
Hence,
for each gauge $q$
on~$E$, there is a gauge $p$
such that
$\|x+y\|_q \leq\max\{ \|x\|_p, \|y\|_p \}$
and thus $\|x+y\|_q \leq \|x\|_p+\|y\|_p$,
for all $x,y\in E$.\,\Punkt
\end{la}
\begin{defn}
Let $E$ be a topological vector space over
a valued field~$\K$.
We say that a set $\Gamma$ of gauges on~$E$ is
a \emph{fundamental system of gauges}
if finite intersections of sets of the form
$B_r^q(0)$ with $r>0$, $q \in\Gamma$
form a basis for the filter of $0$-neighbourhoods
in~$E$.
\end{defn}
Thus, a topological vector space over a valued
field is polynormed if and only if it has a fundamental
system of gauges which are continuous
seminorms. We also mention
that, in the real case,
the continuous gauges always
form a fundamental system
(cf.\ \cite[\S\,6.4]{Jar}).
For a more concrete example,
consider $\ell^p(\K)$ with $p\in \;]0,1]$.
Then $\{\|.\|_p\}$ is a fundamental
system of gauges.\\[2.5mm]
It is useful to know good fundamental
systems of gauges for function spaces.
%
%
\begin{la}\label{spla1}
Let $\K$ be a valued field, $U\sub \K$ be a non-empty subset
without isolated points,
$E$ be a topological $\K$-vector space,
$q$ be a gauge on~$E$,\linebreak
$k\in \N_0\cup\{\infty\}$
and $j\in \N_0$ such that $j\leq k$. Then
%
\begin{equation}\label{standgau}
BC^k(U,E)\to [0,\infty[,\quad
\gamma\mto \|\gamma^{<j>}\|_{q,\infty}:=\,\sup
\{\|\gamma^{<j>}(x)\|_q\colon x\in U^{j+1}\}
\end{equation}
is a gauge on $BC^k(U,E)$.
If $\Gamma$ is a fundamental system
of gauges for~$E$,
then the gauges $\gamma\mto \|\gamma^{<j>}\|_{q,\infty}$
$($for $j\in \N_0$ such that $j\leq k$
and $q\in \Gamma)$
form a fundamental system of gauges
for $BC^k(U,E)$.
\end{la}
If $E=\K$, we simply write $\|\gamma^{<k>}\|_\infty$
instead of $\|\gamma^{<k>}\|_{|.|,\infty}$.
\begin{rem}\label{oldnew2}
If $q$ in Lemma~\ref{spla1} is an upper
semicontinuous gauge which does not happen
to be a seminorm,
then one cannot expect
that the gauge
$\gamma\mto \|\gamma^{<k>}\|_{q,\infty}$
is upper semicontinuous.
For this reason,
we found it convenient to give
up upper semicontinuity in our definition
of gauges. Of course, alternatively
one might redefine
$\|\gamma^{<k>}\|_{q,\infty}$
in a way which enforces
upper semicontinuity,
but such variants
would be more complicated
to work with.
\end{rem}
We need to know
how translations and homotheties affect
the gauges from~(\ref{standgau}).
%
%
\begin{la}\label{transhomo}
Let $\K$ be a valued field,
$U\sub \K$ be a non-empty subset
without isolated points,
$E$ be a topological $\K$-vector space,
$q$ a gauge on~$E$ and $k\in \N_0$.
\begin{itemize}
\item[\rm(a)]
If $\gamma\in BC^k(U,E)$
and $t_0\in\K$, then
$\eta\colon U-t_0\to E$, $\eta(t):=\gamma(t+t_0)$
belongs to $BC^k(U-t_0,E)$.
Furthermore, $\|\eta^{<k>}\|_{q,\infty}=
\|\gamma^{<k>}\|_{q,\infty}$.
\item[\rm(b)]
If $\gamma\in BC^k(U,E)$
and $a\in \K^\times$, then
$\eta\colon a^{-1}U \to E$, $\eta(t):=\gamma(at)$
belongs to $BC^k(a^{-1}U,E)$.
Furthermore,
$\|\eta^{<k>}\|_{q,\infty}=
|a|^k \|\gamma^{<k>}\|_{q,\infty}$.
\item[\rm(c)]
Let $V\sub U$ be a non-empty subset
without isolated points.
Then\linebreak
$\gamma|_V\in BC^k(V,E)$
for $\gamma\in BC^k(U,E)$,
and $\|(\gamma|_V)^{<k>}\|_{q,\infty}\leq
\|\gamma^{<k>}\|_{q,\infty}$.
\end{itemize}
\end{la}
In\hspace*{-.13mm}
\hspace*{-.13mm}the
\hspace*{-.13mm}real
\hspace*{-.13mm}locally
\hspace*{-.13mm}convex
\hspace*{-.13mm}case,
\hspace*{-.2mm}$BC^k$-maps
\hspace*{-.13mm}on
\hspace*{-.13mm}intervals
\hspace*{-.13mm}are
\hspace*{-.13mm}what
\hspace*{-.13mm}they
\hspace*{-.13mm}should~be.
%
%
\begin{la}\label{charctbd}
Let $E$ be a real locally convex space,
$I\sub \R$ be a non-singleton
interval, $k\in \N_0$
and $\gamma\colon I\to E$ be a map.
Then
$\gamma\in C^k(I,E)$
if and only if $\gamma$ is $C^k$
in the usual sense $($viz.\
$\gamma^{(j)}$ exists for $j\in\{0,1,\ldots, k\}$
and is continuous$)$.
Moreover,
$\gamma\in BC^k(I,E)$
if and only if $\gamma$ is $C^k$
in the usual sense
and $\gamma^{(j)}(I)$ is bounded in~$E$
for each $j\in \{0,1,\ldots, k\}$.
\end{la}
\section{Calibrations and absolute convergence}\label{seccalib}
%
%
In this section,
$E$ is a topological vector space
over a valued field~$\K$.
Our goal is to define a meaningful
notion of absolute convergence of
series in~$E$.
As a tool, calibrations
are introduced, which are certain
sequences of gauges.
Compare~\cite{AEK}
for the related concept of a ``string''
in a real vector space.
%
%
\begin{defn}\label{defcalib}
A sequence $(q_n)_{n\in \N_0}$
of gauges on~$E$ is called a \emph{calibration} if
%
\begin{equation}\label{ptycalib}
(\forall n\in\N_0)(\forall x,y\in E)\quad
q_n(x+y)\;\leq\; q_{n+1}(x)+q_{n+1}(y)\,.
\end{equation}
The sequence is a \emph{strong calibration} if
%
\begin{equation}\label{fake2}
(\forall n\in\N_0)(\forall x,y\in E)\quad
q_n(x+y)\;\leq\; \max\{q_{n+1}(x),q_{n+1}(y)\}\,.
\end{equation}
\end{defn}

\noindent We shall refer to (\ref{ptycalib})
as the \emph{fake triangle inequality}.
Similarly, (\ref{fake2}) is called
the \emph{fake ultrametric inequality}.
If $q$ is a gauge on~$E$, then there always
exists a calibration $(q_n)_{n\in \N_0}$
such that $q_0=q$
(cf.\ Lemma~\ref{substitut}).
In this situation, we say that $q$ extends to $(q_n)_{n\in \N_0}$.\\[2.5mm]
In this paper, we decided to work entirely
with ordinary calibrations.
Using strong calibrations instead, one obtains
analogous results. For example,
variants of Lemmas~\ref{estext}
and~\ref{realextzero}
hold for strong calibrations
(in which case the factors $2^{k-j}$
in (\ref{resestim}) and (\ref{morkll}) can be omitted).
%
%
\begin{rem}\label{increas}
If $(q_n)_{n\in \N_0}$ is a  calibration,
then $q_n\leq q_{n+1}$ for each $n\in \N_0$
because $q_n(x)=q_n(x+0)\leq q_{n+1}(x)+q_{n+1}(0)=q_{n+1}(x)$
for each $x\in E$.
\end{rem}
%
%
\begin{rem}\label{excalib}
If $q\colon E\to [0,\infty[$ is a continuous
seminorm, then $(q)_{n\in \N_0}$
is a calibration
(and a strong calibration if $q$ is an ultrametric
seminorm). If $(q_n)_{n\in \N_0}$
is any calibration extending the seminorm~$q$,
then $q_n\geq q$ for each $n$,
by the preceding remark.
Thus $(q)_{n\in \N_0}$ is the smallest
calibration extending~$q$.
\end{rem}
To illustrate the notion of a calibration,
let us look at another example.
\begin{example}
If $r\in \;]0,1]$ and $q$
is an $r$-seminorm on~$E$,
define $q_n:=2^{\frac{n}{r}}q$
for $n\in \N_0$.
Then $(q_n)_{n\in \N_0}$
is a strong calibration on~$E$.
Notably, $(2^{\frac{n}{p}}\|.\|_p)_{n\in \N_0}$
is a strong calibration on $\ell^p(\K)$,
for each $p\in \;]0,1]$.
This follows from the observation
that $q(x+y)^r\leq q(x)^r+q(y)^r\leq 2\max\{q(x)^r,q(y)^r\}$
for $x,y\in E$
and thus $q(x+y)\leq \sqrt[r]{2\max\{q(x)^r,q(y)^r\}}
=\max\{2^{\frac{1}{r}}q(x),2^{\frac{1}{r}}q(y)\}$.
\end{example}
The following lemma is obvious.
\begin{la}\label{spla2}
Let $(q_n)_{n\in \N_0}$
be a calibration on~$E$, $k\in \N_0\cup\{\infty\}$
and \mbox{$j\in \N_0$}
with $j\leq k$.
Then the gauges $BC^k(U,E)\to[0,\infty[$,
$\gamma\mto \|\gamma^{<j>}\|_{q_n,\infty}$,
for $n\in \N_0$,
form a calibration on $BC^k(U,E)$.\Punkt
\end{la}
Calibrations are valuable tools to establish
the convergence of series in topological
vector spaces which may fail to be polynormed.
%
\begin{defn}\label{defabscon}
Let $(x_n)_{n\in \N}$
be a sequence in $E$.
We say that the series $\sum_{n=1}^\infty x_n$
is \emph{absolutely convergent} if each
gauge $q$ on~$E$
extends to a calibration $(q_n)_{n\in \N_0}$
such that
\[
\sum_{n=1}^\infty \|x_n\|_{q_n} \; < \; \infty\, .
\]
\end{defn}
\begin{rem}
If~$E$ is polynormed, then a series $\sum_{n=1}^\infty x_n$
in~$E$ is absolutely convergent if and only if
$\sum_{n=1}^\infty \|x_n\|_q<\infty$
for each continuous seminorm~$q$ on~$E$
(cf.\ Remark~\ref{excalib}).
\end{rem}
Absolute convergence of series in a topological vector space
is a useful concept
provided that
the latter
is \emph{sequentially complete}
in the sense that each Cauchy sequence converges.
%
\begin{la}\label{absconcon}
If $E$ is a sequentially complete
topological vector space over a valued field~$\K$,
then every absolutely convergent series
in~$E$
is convergent.
\end{la}
\begin{proof}
Using (\ref{ptycalib}) repeatedly, we see that
$\big\|\sum_{k=m}^n x_k\big\|_{q_0}\leq
\sum_{k=m}^n \|x_k\|_{q_{k-m+1}}$\linebreak
$\leq
\sum_{k=m}^n \|x_k\|_{q_k}$
for all $n ,m\in \N$ with
$n>m$.
This entails that $(\sum_{k=1}^n x_k)_{n\in \N}$
is a Cauchy sequence in~$E$ and hence convergent.
\end{proof}
\section{Ultrametric General Curve Lemma}\label{secgencu}
In this section, we formulate and prove our
first main result.
%
%
\begin{thm}[Ultrametric General Curve Lemma]\label{thmgencu}
Let $E$ be a topological vector space
over an ultrametric field~$\K$,
$\rho\in \K^\times$ with $|\rho|<1$
and $(\gamma_n)_{n\in \N}$ be a sequence
of smooth maps $\gamma_n\in BC^\infty(\rho^n\bO, E)$
which become small sufficiently fast
in the sense that,
for each gauge $q$ on~$E$,
there exists a calibration $(q_n)_{n\in\N_0}$
extending $q$ such that
%
\begin{equation}\label{condgen}
(\forall a>0)\, (\forall k, m\in \N_0)\quad
\lim_{n\to\infty}\,
a^n\|\gamma_n^{<k>}\|_{q_{n+m}, \infty}\;=\;  0\, .
\end{equation}
Then there exists a smooth map
$\gamma\in BC^\infty(\K,E)$
whose image
$\im(\gamma)$ is\linebreak
contained in
$\{0\}\cup\bigcup_{n\in \N}\im(\gamma_n)$,
such that
\begin{equation}\label{desire}
\gamma(\rho^{n-1}+t)\;=\;
\gamma_n(t)\quad\mbox{for all $n\in \N$ and $t\in \rho^n\bO$.}
\end{equation}
\end{thm}
\begin{rem}
Note that $\rho^n\bO= \wb{B}_{|\rho|^n}(0)$
and $\rho^{n-1}+\rho^n\bO=\wb{B}_{|\rho|^n}(\rho^{n-1})$
here.
Since $|\rho|^n<|\rho^{n-1}|$, we have
$|x|=|\rho^{n-1}|$
for each $x\in \wb{B}_{|\rho|^n}(\rho^{n-1})$
(see (\ref{winner})).
As a consequence, the balls
$\wb{B}_{|\rho|^n}(\rho^{n-1})$
are pairwise disjoint.
\end{rem}
\begin{rem}
If~$E$ is polynormed, then
the somewhat complicated condition (\ref{condgen})
can be simplified.
In view of Remark~\ref{excalib},
condition (\ref{condgen}) then amounts
to the following:
For each $k\in \N_0$ and continuous seminorm
$q$ on~$E$, we have
%
\begin{equation}\label{condgen2}
(\forall a>0)\quad\;\,
\lim_{n\to\infty}\, a^n\|\gamma_n^{<k>}\|_{q, \infty}\; = \; 0\,.
\end{equation}
\end{rem}
%
%
\begin{rem}\label{shallrm}
Let $E$ in Lemma~\ref{thmgencu} be metrizable
and suppose that there exists
a calibration $(p_n)_{n\in \N_0}$
such that $\{p_n\colon n\in \N_0\}$
is a fundamental system of gauges,
and $C>0$ such that
%
\begin{equation}\label{shalluse}
(\forall k\in \N_0)\, (\forall n\geq k)\quad
\|\gamma_n^{<k>}\|_{p_{2n},\infty}\;\leq\; C n^{-n}\,.
\end{equation}
Then the hypothesis (\ref{condgen})
of Theorem~\ref{thmgencu} is satisfied:
Given~$q$, we can extend it to a suitable calibration
via $q_n:=rp_{n+n_0}$ for $n \in \N$,
with $r>0$ and $n_0\in \N_0$ sufficiently large.
In all our applications,
we use this simpler~criterion.
\end{rem}
The following lemma prepares the proof of
Theorem~\ref{thmgencu}. As before, $\K$ is an
ultrametric field and~$E$ a topological
$\K$-vector space.
%
%
\begin{la}\label{estext}
Let $\gamma\in BC^\infty(U,E)$,
where $U:=\wb{B}_r^\K(0)$ for some
$r \in \;]0,\infty[$.
Extend $\gamma$ to a smooth map
$\eta\colon \K\to E$ via
$\eta(x):=0$ for $x\in \K\setminus U$.
Then $\eta\in BC^\infty(\K,E)$,
and\vspace{-2.5mm}
\begin{equation}\label{resestim}
\|\eta^{<k>}\|_{q_0,\infty}\;\leq\;
\, \max_{j=0,\ldots, k}\,
\Big(\frac{2}{r}\Big)^{k-j}\,
\|\gamma^{<j>}\|_{q_{k-j},\infty}\,,\vspace{-2mm}
\end{equation}
for each $k\in \N_0$ and calibration $(q_n)_{n\in \N_0}$ on~$E$.
\end{la}
\begin{proof}
Note first that~$\eta$ is smooth since
smoothness is a local property
(see \cite[Lemma~4.9]{BGN})
and $U$ is both open and closed.
We now show by induction on $k\in \N_0$
that $\eta\in BC^k(\K,E)$
and~(\ref{resestim}) holds.
If $k=0$, then
$\eta\in BC(\K,E)$
and (\ref{resestim}) holds
because $\sup\{\|\eta(x)\|_{q_0}\colon x\in \K\}
=\sup\{\|\gamma(x)\|_{q_0}\colon x\in U \}
=\|\gamma^{<k>}\|_{q_0,\infty}$.
Now suppose that $k\geq 1$ and suppose
that (\ref{resestim})
holds if $k$ is replaced with $k-1$,
for each calibration.
Since $U^{>k<}$ is dense in $U^{<k>}$
and $\eta^{<k>}$ is continuous,
we only need to show that
the right hand
side of (\ref{resestim})
is an upper bound for $\|\eta^{<k>}(x)\|_{q_0}$,
for each $x=(x_0,\ldots, x_k)\in U^{>k<}$.
It is convenient to distinguish three cases:\\[2.5mm]
\emph{Case}~1: If $x_j\in U$ for all
$j\in\{0,\ldots,k\}$, then
\[
\eta^{<k>}(x_0,\ldots,x_k)\; =\; \gamma^{<k>}(x_0,\ldots, x_k)
\]
and thus
$\|\eta^{<k>}(x_0,\ldots,x_k)\|_{q_0}
=\|\gamma^{<k>}(x_0,\ldots, x_k)\|_{q_0}
\leq\|\gamma^{<k>}\|_{q_0,\infty}$,
which does not exceed the right hand
side of (\ref{resestim}).\\[2.5mm]
\emph{Case}~2: If $x_j\not\in U$ for all
$j\in\{0,\ldots,k\}$, then
$\eta(x_j)=0$ for each $j$ and thus
$\eta^{<k>}(x_0,\ldots,x_k)=0$,
whence again
$\|\eta^{<k>}(x_0,\ldots,x_k)\|_{q_0}=0$
does not exceed the right hand
side of (\ref{resestim}).\\[2.5mm]
\emph{Case}~3: There are $i,j\in \{0,\ldots, k\}$
such that $x_i\in U$ and $x_j\not\in U$.
By symmetry of $\eta^{<k>}$, we may assume
that $i=0$ and $j=k$.
Since $|x_0|\leq r< |x_k|$
and $|.|$ is ultrametric, we have $|x_0-x_k|=|x_k|>r$.
Hence
\begin{eqnarray}
\|\eta^{<k>}(x)\|_{q_0} &=&
\frac{\|\eta^{<k-1>}(x_0,x_1,\ldots, x_{k-1})-
\eta^{<k-1>}(x_k,x_1,\ldots, x_{k-1})\|_{q_0}}{|x_0-x_k|}
\notag \\
&\leq&
\frac{\|\eta^{<k-1>}(x_0,x_1,\ldots, x_{k-1})\|_{q_1}+
\|\eta^{<k-1>}(x_k,x_1,\ldots, x_{k-1})\|_{q_1}}{r}
\notag \\
&\leq&
\frac{2}{r}\cdot
\max_{j=0,\ldots, k-1}\,
\Big(\frac{2}{r}\Big)^{k-1-j}\|\gamma^{<j>}\|_{q_{k-j},\infty}\,, \label{inddd}
\end{eqnarray}
applying
the inductive hypothesis
to $\eta^{<k-1>}$ and the calibration
$(q_{n+1})_{n\in \N_0}$
to obtain the final inequality.
Since the right hand side of (\ref{inddd})
does not exceed the right hand
side of (\ref{resestim}),
our inductive proof is complete.
\end{proof}
{\bf Proof of Theorem~\ref{thmgencu}.}
For each $n\in \N$, define
$\eta_n\colon \K\to E$ via
\[
\eta_n(t)\;:=\;
\left\{
\begin{array}{cl}
\gamma_n(t-\rho^{n-1}) & \;\mbox{if $\,|t-\rho^{n-1}|\leq\rho^n$;}\\
0 & \; \mbox{otherwise.}
\end{array}
\right.
\]
Then $\eta_n\in BC^\infty(\K,E)$
and
%
\begin{equation}\label{uglyest}
\|\eta_n^{<k>}\|_{q_0,\infty} \; \leq \;
\max_{j=0,\ldots, k}\,
\big({\textstyle \frac{2}{|\rho|^n}}\big)^{k-j}
\|\gamma_n^{<j>}\|_{q_{k-j},\infty}
\end{equation}
for each $k\in \N_0$ and calibration $(q_n)_{n\in \N_0}$
on~$E$, by Lemmas~\ref{estext}
and~\ref{transhomo}\,(a).
Define $\gamma(x):=\sum_{n=1}^\infty\eta_n(x)$
for $n\in \N$.
Then $\gamma\colon \K\to E$ is smooth
on $\K^\times$, using that the supports
of the maps $\eta_n$ form a locally finite
family of disjoint subsets of $\K^\times$
(since
$\Supp(\eta_n)\sub \wb{B}_{|\rho|^n}(\rho^{n-1})$).

\emph{Step}~1: We show that
$\sum_{n=1}^\infty\eta_n$
converges in $BC^\infty(\K,\wb{E})$,
where $\wb{E}$ is the completion of~$E$.
Once this is established,
for each $x\in U$ we
can apply\linebreak
the continuous linear point evaluation
$BC^\infty(\K,\wb{E})\to \wb{E}$,
$\zeta\mto\zeta(x)$
to $\sum_{n=1}^\infty\eta_n$,
showing that $(\sum_{n=1}^\infty\eta_n)(x)=\gamma(x)$.
Since $BC^\infty(\K,\wb{E})$
is complete by Lemma~\ref{completeBC}\,(d),
to establish convergence we only need to show
that the series $\sum_{n=1}^\infty\eta_n$
converges absolutely
(see Lemma~\ref{absconcon}).
To this end, let $q_0$ be a gauge on~$\wb{E}$
and extend it to
a calibration
$(q_n)_{n\in \N_0}$ such that (\ref{condgen})
holds (and hence also (\ref{uglyest})).
Since the gauges $BC^\infty(\K,\wb{E})\to [0,\infty[$,
$\zeta\mapsto \|\zeta^{<k>}\|_{q,\infty}$
form a fundamental system of gauges for
$k$ ranging through $\N_0$ and $q$ through the gauges
of~$\wb{E}$ (see Lemma~\ref{spla1}),
$\sum_{n=1}^\infty\eta_n$
will converge absolutely in $BC^\infty(\K,\wb{E})$
if we can show
that $\sum_{n=1}^\infty \|\eta_n^{<k>}\|_{q_n,\infty}<\infty$
in the preceding situation,
for each $k\in \N_0$.
In view of (\ref{uglyest}), it suffices to show that
\[
\sum_{n=1}^\infty \max_{j=1,\ldots, k}
\big({\textstyle \frac{2}{|\rho|^n}}\big)^{k-j}
\|\gamma_n^{<j>}\|_{q_{n+k-j},\infty}\; <\;\infty\,.
\]
This will hold if we can show that, for each $j\in
\{1,\ldots, k\}$,
\begin{equation}\label{needshow}
\sum_{n=1}^\infty|\rho|^{-nk}\|\gamma_n^{<j>}\|_{q_{n+k-j},\infty}
\;<\;\infty\,.
\end{equation}
To prove (\ref{needshow}), choose
$a>|\rho|^{-k}$ and recall that
$a^n\|\gamma_n^{<j>}\|_{q_{n+k-j},\infty}\to 0$,
by~(\ref{condgen}).
Thus $A_k:=\sup\,\{a^n\|\gamma_n^{<j>}\|_{q_{n+k-j},\infty}\colon n\in \N\}
<\infty$
and hence $\sum_{n=1}^\infty A_k(\frac{|\rho|^{-k}}{a})^n$
is a convergent majorant for
$\sum_{n=1}^\infty|\rho|^{-nk}
\|\gamma_n^{<j>}\|_{q_{n+k-j},\infty}$.\\[2.5mm]
\emph{Step}~2: We now show by induction
on $k\in \N_0$ that $\gamma^{<k>}\colon \K^{k+1}\to \wb{E}$
actually takes values in~$E$,
whence $\gamma$ is $C^k$ as a map into~$E$.
For $k=0$, this is trivial.
Let $k\geq 1$ now and assume that
the assertion holds if $k$ is replaced with $k-1$.
Let $x=(x_0,\ldots, x_k)\in \K^{k+1}$.
If $x_i\not=x_j$ for certain $i,j\in \{0,\ldots, k\}$,
then $\gamma^{<k>}(x)$ is a partial difference
quotient of $\gamma^{<k-1>}$ and hence a scalar
multiple of two values of $\gamma^{<k-1>}$,
which are in~$E$ (by induction).
Hence also $\gamma^{<k>}(x)\in E$.
It remains to show that
$\gamma^{<k>}(y,\ldots, y)\in E$
for all $y\in \K$. If $y\in \K^\times$,
this follows from the smoothness of $\gamma|_{\K^\times}$. 
To see that $\gamma^{<k>}(0)\in E$,
we exploit the continuity of the map
\[
BC^\infty(\K,\wb{E})\to E\, , \quad
\zeta\mapsto \zeta^{<k>}(0,\ldots, 0)\,.
\]
It entails that
$\gamma^{<k>}(0)=\sum_{n=1}^\infty
\eta_n^{<k>}(0)=0\in E$.
Thus $\gamma^{<k>}$ has image in~$E$,
which completes the induction.\vspace{3mm}\Punkt

\noindent
%
%
%
%
%
%
%
%
%
%
%
%
%
\section{General
Curve Lemma for curves in
real topological vector spaces}\label{gencureal}
%
%
%
In this section, we prove a version of
the General Curve Lemma for curves in arbitrary
(not necessarily locally convex) real
topological vector spaces.
%
%
\begin{thm}[Real Case of General Curve Lemma]\label{thmrealcu}
Let $E$ be a real topological
vector space and $(s_n)_{n\in \N}$
as well as $(r_n)_{n\in \N}$
be sequences of positive reals such that
$\sum_{n=1}^\infty s_n<\infty$
and $r_n\geq s_n+\frac{2}{n^2}$
for each $n\in \N$.
Furthermore, let
$(\gamma_n)_{n\in \N}$ be a sequence of smooth
maps $\gamma_n\colon [-r_n,r_n]\to E$
which become small sufficiently
fast in the sense that, for each gauge~$q$ on~$E$,
there exists a calibration $(q_n)_{n\in \N_0}$
extending $q$ such that
%
\begin{equation}\label{dfco}
(\forall k,\ell, m\in \N_0)\qquad
\lim_{n\to\infty} \, n^\ell\|\gamma_n^{<k>}\|_{q_{n+m},\infty}\,=\,0\,.
\end{equation}
Then there exists a curve $\gamma\in BC^\infty(\R,E)$
with $\,\im(\gamma)\sub [0,1]\cdot\bigcup_{n\in \N}\im(\gamma_n)$
and a convergent sequence $(t_n)_{n\in \N}$
of real numbers such that
\begin{equation}\label{nicereal}
\gamma(t_n+t)\;=\; \gamma_n(t)\quad
\mbox{for all $\, n\in \N$ and $\,t\in [-s_n,s_n]$.}
\end{equation}
\end{thm}
Various lemmas are needed to prepare the proof
of Theorem~\ref{thmrealcu}.
%
%
\begin{la}\label{enespro}
For each $n\in \N$,
there exist integers $N_{i,j}\in \N_0$
indexed by all strictly increasing finite
sequences
$i=(i_0,\ldots,i_k)$
and $j=(j_0,\ldots, j_\ell)$
with entries in $\{0,1,\ldots, n\}$,
for $k,\ell\in \N_0$ with $k+\ell=n$,
such that $\sum_{i,j}N_{i,j}\leq 2^n$
and the following holds:
For each topological field~$\K$,
non-empty subset $U\sub\K$ without isolated
points, continuous bilinear map
$\beta\colon E\times F\to H$
between topological $\K$-vector spaces
and all $C^n$-maps
$\gamma\colon U\to E$,
$\eta\colon U\to F$,
we have
%
\begin{eqnarray}
\hspace*{-10mm}\lefteqn{\big(\beta\circ (\gamma,\eta)\big)^{<n>}
(x_0,\ldots, x_n)}\qquad\notag\\
&= &
\sum_{k+\ell=n}\sum_{\stackrel{{\scriptstyle i,j\;\text{with}}}{\#i=k,\;\#j=\ell}}
N_{i,j}\;
\beta\big(
\gamma^{<k>}(x_{i_0},\ldots, x_{i_k}),\,
\eta^{<\ell>}(x_{j_0},\ldots, x_{j_\ell})\big)\label{goodenough}
\end{eqnarray}
for all $(x_0,\ldots, x_n)\in U^{n+1}$,
using the notation $\# (i_0,\ldots, i_k):=k$.
\end{la}
%
%
\begin{rem}\label{meani}
The condition
$\sum_{i,j}N_{i,j}\leq 2^n$
means that we can consider\linebreak
$\big(\beta\circ (\gamma,\eta)\big)^{<n>}
(x_0,\ldots, x_n)$
as a sum of $\leq 2^n$ summands
of the form
\[
\beta\big(
\gamma^{<k>}(x_{i_0},\ldots, x_{i_k}),\,
\eta^{<\ell>}(x_{j_0},\ldots, x_{j_\ell})\big)\,.
\]
\end{rem}
{\bf Proof of Lemma~\ref{enespro}.}
The proof is by induction on~$n\in \N$.
If $n=1$ and $x_0,x_1\in U$ are distinct, then
%
\begin{eqnarray}
\lefteqn{\frac{\beta(\gamma(x_1),\eta(x_1))
-\beta(\gamma(x_0),\eta(x_0))}{x_1-x_0}}\qquad\notag \\
&=&
\frac{\beta(\gamma(x_1),\eta(x_1))-\beta(\gamma(x_0),\eta(x_1))+
\beta(\gamma(x_0),\eta(x_1))-
\beta(\gamma(x_0),\eta(x_0))}{x_1-x_0}\notag \\
&=&
\beta(\gamma^{<1>}(x_0,x_1),\eta(x_1))
+\beta(\gamma(x_0),\eta^{<1>}(x_0,x_1))\,.\label{extable}
\end{eqnarray}
Since (\ref{extable}) can be used to define
a continuous function in $(x_0,x_1)\in U\times U$,
we see that $\beta\circ (\gamma,\eta)$
is $C^1$ with
%
\begin{equation}\label{csn1}
(\beta\circ(\gamma,\eta))^{<1>}(x_0,x_1)=
\beta(\gamma^{<1>}(x_0,x_1),\eta(x_1))
+\beta(\gamma(x_0),\eta^{<1>}(x_0,x_1))
\end{equation}
of the form described in (\ref{goodenough}).\\[3mm]
\emph{Induction step}: Suppose that the lemma
holds for some~$n$ and that $\gamma,\eta$ are $C^{n+1}$.
For $i,j$ as above with $\# i=k$,
$\# j =\ell$ and $k+\ell=n$,
abbreviate
\[
h_{i,j}(x_0,\ldots, x_n)\,:=\,
\beta\big(
\gamma^{<k>}(x_{i_0},\ldots, x_{i_k}),\,
\eta^{<\ell>}(x_{j_0},\ldots, x_{j_\ell})\big)
\]
for $x_0,\ldots, x_n\in U$.
The analogue of (\ref{goodenough})
for $(\beta\circ(\gamma,\eta))^{<n+1>}$
will be apparent from an explicit
formula for the continuous extension
of the mapping\linebreak
$g\colon \{x=(x_0,\ldots, x_{n+1})\in U^{n+2}\colon x_0\not=x_{n+1}\}
\to H$,
%
\begin{equation}\label{prlm}
g(x)\;:=\;
\frac{h_{i,j}(x_{n+1},x_1,\ldots, x_n)-h_{i,j}(x_0,x_1,\ldots,
x_n)}{x_{n+1}-x_0}
\end{equation}
to a map $U^{n+2}\to H$,
which we now establish.
If $i_0\not=0$ and $j_0\not=0$,
then $h_{i,j}$ does not depend on $x_0$
and thus $g=0$ has $0$ as a continuous extension.\\[3mm]
If $i_0=0$ and $j_0\not=0$,
then $\eta^{<\ell>}(x_{j_0},\ldots, x_{j_\ell})$
does not depend on~$x_0$ and thus
\[
g(x_0,\ldots,x_{n+1})
\; =\;
\beta(\gamma^{<k+1>}(x_{i_0},\ldots, x_{i_k},x_{n+1}),
\eta^{<\ell>}(x_{j_0},\ldots, x_{j_\ell}))
\]
by linearity of~$\beta$ in its first argument,
where the right hand side
can be used to define a continuous function
on $U^{n+2}$.
Likewise, the mapping
$U^{n+2}\to H$,
$x\mto \beta(\gamma^{<k>}(x_{i_0},\ldots, x_{i_k}),
\eta^{<\ell+1>}(x_{j_0},\ldots, x_{j_\ell},x_{n+1}))$
provides a continuous extension of~$g$
if $i_0\not=0$ and $j_0=0$.\\[3mm]
If $i_0=j_0=0$, then
the calculation leading to (\ref{extable})
shows that
\begin{eqnarray*}
g(x) &= & \beta(\gamma^{<k+1>}(x_{i_0},\ldots,x_{i_k},x_{n+1}),
\eta^{<\ell>}(x_{j_1},\ldots,x_{j_\ell},x_{n+1}))\\
& & \;\; +
\beta(\gamma^{<k>}(x_{i_0},\ldots, x_{i_k}),
\eta^{<\ell+1>}(x_{j_0},\ldots, x_{j_\ell},x_{n+1}))\,,
\end{eqnarray*}
where again the right hand side
extends continuously
to all of $U^{n+2}$.\\[3mm]
Forming the sum of all contributions
just described,  we obtain
a formula analogous to (\ref{goodenough})
for $(\beta\circ(\gamma,\eta))^{<n+1>}$.\vspace{2.5mm}\Punkt

\begin{rem}\label{ludkuexpl}
For our purposes, we need not know
the integers $N_{i,j}$ explicitly.
\end{rem}
%
%
%
\begin{la}\label{realestim}
There exist
constants $C_k \in \N$
for $k\in \N_0$
with $\sum_{k=0}^nC_k\leq 2^n$
for each $n\in \N$,
and the following property:
For each valued field $\K$,
non-empty subset $U\sub \K$ without isolated
points, $n\in \N$, $\gamma\in BC^n(U,\K)$,
topological $\K$-vector space~$E$,
$\eta\in BC^n(U,E)$
and calibration $(q_k)_{k\in \N_0}$ on~$E$, we have
%
\begin{equation}\label{estprdreal}
\|(\gamma\cdot \eta)^{<n>}\|_{q_0,\infty}
\;\leq\;\, \sum_{k=0}^n \,C_k\,\|\gamma^{<k>}\|_\infty
\cdot \|\eta^{<n-k>}\|_{q_n,\infty}\,.
\end{equation}
\end{la}
\begin{proof}
Applying (\ref{goodenough})
to the scalar multiplication $\beta\colon \K\times E\to E$,
we get
\[
\|(\gamma\cdot \eta)^{<n>}\|_{q_0,\infty}\;\leq\;
\sum_{k=0}^n\sum_{i,j} \, |N_{i,j}|\hspace*{.2mm}
\|\gamma^{<k>}\|_\infty
\|\eta^{<n-k>}\|_{q_n,\infty}\,.
\]
Here, $\leq 2^n$ summands were involved and
hence the fake triangle inequality had to be used
at most $n$ times, explaining why the
gauge $q_n$ occurs.
Since $|N_{i,j}|\leq N_{i,j}$,
the assertion follows with
$C_k:=\sum_{i,j}N_{i,j}$,
where the sum is taken over all $i,j$
as in Lemma~\ref{enespro} such that $\# i=k$
and $\# j=n-k$.
\end{proof}
As a first application of Lemma~\ref{realestim},
let us construct a family of smooth cut-off functions
the size of whose
difference quotient maps
(of all orders) is well under control.
These cut-off functions will be most
useful later.
%
%
\begin{la}\label{scutoff}
There is a sequence $(M_n)_{n\in \N_0}$
of positive reals with the following property:
For all $a,b>0$, there exists a smooth function
$h\colon \R\to[0,1]$
with support $\Supp(h)\sub [{-(a+b)},a+b]$,
such that
$h(t)=1$ for all $t\in [{-a},a]$ and
\[
\hspace*{-13mm}(\forall n\in \N_0)\qquad
\|h^{<n>}\|_\infty \;\leq\; M_n\, b^{-n}\,.
\]
\end{la}
\begin{proof}
Let $g\colon \R\to [0,1]$ be a smooth function
such that $g(t)=1$ if $t\leq 0$ and $g(t)=0$ if $t\geq 1$. 
Then $g^{(k)}$ is bounded for each $k\in \N_0$
and hence $g\in BC^\infty(\R,\R)$, by Lemma~\ref{charctbd}.
Given $a,b>0$, define $h\colon \R\to\R$
via $h(t):=g(\frac{t-a}{b})g(\frac{-t-a}{b})$.
Then $h(\R)\sub [0,1]$,
$h(t)=1$ if $|t|\leq a$, and $h(t)=0$
if $|t|\geq a+b$.
By Lemma~\ref{realestim},
we have $h\in BC^\infty(\R,\R)$.
Furthermore, combining (\ref{estprdreal}) with
Lemma~\ref{transhomo} (a) and (b), we see that\vspace{-2mm}
\[
\|h^{<n>}\|_\infty
\;\leq\;\, \sum_{k=0}^n \,C_k\,b^{-k}
\|g^{<k>}\|_\infty \cdot b^{-(n-k)}\|g^{<n-k>}\|_\infty
\;=\; M_n\, b^{-n}\vspace{-1.5mm}
\]
with
$\,M_n:= \sum_{k=0}^n C_k \|g^{<k>}\|_\infty \cdot
\|g^{<n-k>}\|_\infty$ independent of $a$ and~$b$.
\end{proof}
The following lemma will serve as a substitute
for Lemma~\ref{estext} in the real case.
Of course,
$BC^\infty([a,b], E)
=C^\infty([a,b], E)$
here by compactness of~$[a,b]$.
%
%
\begin{la}\label{realextzero}
Let $a<\alpha<\beta<b$ be real numbers,
$r:=\min\{\alpha-a, b-\beta\}$,
$E$ be a real topological vector space,
and $\gamma\in BC^\infty([a,b], E)$
be a map such that $\gamma(x)=0$ if $x\in [a,b]\setminus [\alpha,\beta]$.
Define $\eta\colon \R\to E$ via $\eta(x):=\gamma(x)$
if $x\in [a,b]$, $\eta(x):=0$ else.
Then $\eta\in BC^\infty(\R,E)$.
Furthermore,\vspace{-1.5mm}
%
\begin{equation}\label{morkll}
\|\eta^{<k>}\|_{q_0,\infty}\;\leq\;
\max_{j=0,\ldots, k}\,
\Big(\frac{2}{r}\Big)^{k-j}
\|\gamma^{<j>}\|_{q_{k-j},\infty}\, ,\vspace{-1mm}
\end{equation}
for each $k\in \N_0$ and
calibration $(q_n)_{n\in \N_0}$ on~$E$.
\end{la}
\begin{proof}
We show by induction on $k\in \N_0$
that $\eta\in BC^k(\R,E)$
and (\ref{morkll}) holds.
If $k=0$, then
$\sup\{\|\eta(x)\|_{q_0}\colon x\in \R\}
=\sup\{\|\gamma(x)\|_{q_0}\colon x\in [a,b]\}
=\|\gamma^{<k>}\|_{q_0,\infty}$
for each calibration $(q_n)_{n\in \N_0}$,
entailing that $\eta\in BC(\R,E)$ and (\ref{morkll}) holds.
Now suppose that $k\geq 1$ and suppose
that the estimate (\ref{morkll})
holds if $k$ is replaced with $k-1$,
for each calibration.
Since $U^{>k<}$ is dense in $U^{<k>}$
and $\eta^{<k>}$ is continuous,
we only need to show that
the right hand
side of (\ref{morkll})
is an upper bound for $\|\eta^{<k>}(x)\|_{q_0}$,
for each $x=(x_0,\ldots, x_k)\in U^{>k<}$.\\[2.5mm]
\emph{Cases}~1 and~2:
If $x_j\in [a,b]$ for all
$j\in\{0,\ldots,k\}$,
or if $x_j\not\in [\alpha,\beta]$ for all
$j\in\{0,\ldots,k\}$, then
we see as in Step~1 and~2
of the proof of Lemma~\ref{estext}
that $\|\eta^{<k>}(x)\|_{q_0}$
does not exceed the right hand
side of (\ref{morkll}).\\[2.5mm]
\emph{Case}~3: Assume that there are
$i,j\in \{0,\ldots, k\}$
such that $x_i\in [\alpha,\beta]$ and $x_j\not\in [a,b]$.
Then $|x_i-x_j|\geq r$.
By symmetry of $\eta^{<k>}$, without loss of generality
$i=0$ and $j=k$. Now
\begin{eqnarray}
\|\eta^{<k>}(x)\|_{q_0} &=&
\frac{\|\eta^{<k-1>}(x_0,x_1,\ldots, x_{k-1})-
\eta^{<k-1>}(x_k,x_1,\ldots, x_{k-1})\|_{q_0}}{|x_0-x_k|}
\notag \\
&\leq&
\frac{\|\eta^{<k-1>}(x_0,x_1,\ldots, x_{k-1})\|_{q_1}+
\|\eta^{<k-1>}(x_k,x_1,\ldots, x_{k-1})\|_{q_1}}{r}
\notag \\
&\leq&
\frac{2}{r}\cdot
\max_{j=0,\ldots, k-1}\,
\Big(\frac{2}{r}\Big)^{k-1-j}\|\gamma^{<j>}\|_{q_{k-j},\infty}\,, \label{inddd2}
\end{eqnarray}
applying
the inductive hypothesis
to $\eta^{<k-1>}$ and the calibration
$(q_{n+1})_{n\in \N_0}$
to obtain the final inequality.
Since the right hand side of (\ref{inddd2})
does not exceed the right hand
side of (\ref{morkll}),
our inductive proof is complete.
\end{proof}
{\bf Proof of Theorem~\ref{thmrealcu}.}
After shrinking $r_n$ if necessary, we may assume that
$r_n=s_n+\frac{2}{n^2}$ for each $n\in \N$
(cf.\ Lemma~\ref{transhomo}\,(c)).
Let $(M_n)_{n\in \N_0}$ be as in
Lemma~\ref{scutoff}.
Given $n\in \N$,
we apply Lemma~\ref{scutoff} with $a:=s_n$ and $b:=\frac{1}{n^2}$.
We obtain a smooth function
$h_n\colon \R\to [0,1]$
such that $h_n(t)=1$ for all $t\in [{-s_n},s_n]$,
$\Supp(h_n)\sub [{-s_n-\frac{1}{n^2}},s_n+\frac{1}{n^2}]$, and
%
\begin{equation}\label{forl}
\|h_n^{<k>}\|_\infty\;\leq \; M_k n^{2k}
\quad\mbox{for each $k\in \N_0$.}
\end{equation}
Set $r_0:=0$ and define for $n\in \N$\vspace{-2.5mm}
\[
t_n\; :=\; \sum_{j=1}^n (r_j+r_{j-1})\,.\vspace{-2mm}
\]
Then $(t_n)_{n\in\N}$ is a monotonically increasing
sequence, which converges
because
$t_\infty:=\sum_{j=1}^\infty(r_j+r_{j-1})\leq
2\sum_{j=1}^\infty s_j+ 4\sum_{j=1}^\infty\frac{1}{n^2}
<\infty$.
By definition,
%
\begin{equation}\label{givesdisj}
(\forall n\in \N)\qquad t_{n+1}-t_n\,= \, r_{n+1}+r_n\,.
\end{equation}
Define $\zeta_n\colon [t_n-r_n,t_n+r_n]\to E$,
$\zeta_n(t):=h_n(t-t_n)\gamma_n(t-t_n)$
and let $\eta_n\colon \R\to E$ be the extension of~$\zeta_n$
by~$0$.
Then
$\Supp(\eta_n)
\sub [t_n-s_n-\frac{1}{n^2},t_n+s_n+
\frac{1}{n^2}]\sub$\linebreak
$]t_n-r_n,t_n+r_n[$,
whence
the maps $\eta_n$ have disjoint supports
(cf.\ (\ref{givesdisj})).
Thus $\gamma(t):=\sum_{n=1}^\infty\eta_n(t)$
exists pointwise. To see that $\gamma$
has the desired properties,
let $q$ be a gauge on~$E$ and extend
it to a calibration~$(q_n)_{n\in \N_0}$ such that (\ref{dfco})
holds. Then\vspace{-3.5mm}
%
\begin{eqnarray}
\|\zeta_n^{<j>} \|_{q_m,\infty} &\leq &
\sum_{i=0}^j C_i\,\|h_n^{<i>}\|_\infty\cdot
\|\gamma_n^{<j-i>} \|_{q_{m+j},\infty}\notag \\
&\leq&
\sum_{i=0}^j n^{2i} C_iM_i 
\|\gamma_n^{<j-i>}\|_{q_{m+j},\infty}\, ,\vspace{-1mm} \label{maybeus}
\end{eqnarray}
for all $n\in \N$ and $m,j \in \N_0$,
using Lemma~\ref{transhomo}\,(a), inequality~(\ref{estprdreal})
from Lemma~\ref{realestim} and (\ref{forl}).
Since $\zeta_n$ vanishes outside
$[t_n-s_n-\frac{1}{n^2}, t_n+s_n+\frac{1}{n^2}]$
and furthermore
$(t_n-s_n-\frac{1}{n^2})-(t_n-r_n)=
\frac{1}{n^2}$
and $(t_n+r_n)-(t_n+s_n+\frac{1}{n^2})=\frac{1}{n^2}$,
Lemma~\ref{realextzero} and (\ref{morkll})
show that $\eta_n\in BC^\infty(\R,E)$, with\vspace{-2mm}
\begin{eqnarray*}
\|\eta_n^{<k>}\|_{q_n,\infty} \!\! &\! \leq \!\! &\!\!
\max_{j=0,\ldots, k}\, \Big(\frac{2}{1/n^2}\Big)^{k-j}
\|\zeta_n^{<j>}\|_{q_{n+k-j},\infty}
\leq
2^kn^{2k}
\max_{j=0,\ldots, k} \,
\|\zeta_n^{<j>}\|_{q_{n+k},\infty}\\
\!\! &\!\! \leq \!\! & \!
n^{4k}A_k{\textstyle \sum_{i=0}^k}\, \|\gamma_n^{<i>}\|_{q_{n+2k},\infty}
\end{eqnarray*}
where $A_k:=\max\{2^kC_jM_j \colon j\!=\!0,\ldots, k\}$.
Passing to the last line,
we~used (\ref{maybeus})
and replaced some terms by larger ones.
Since $n^{4k+2}\|\gamma_n^{<i>}\|_{q_{n+2k},\infty}$
converges as $n\to\infty$
for each $i\in \{0,\ldots, k\}$
(by~(\ref{dfco})), we have
$B_k := \sup\{n^{4k+2}
\|\gamma_n^{<i>}\|_{q_{n+2k},\infty}\colon
\mbox{$i\in\{0,\ldots, k\}$,
$n\in \N$}\}< \infty$.
Hence\vspace{-2.4mm}
\[
\sum_{n=1}^\infty \|\eta_n^{<k>}\|_{q_n,\infty}
\leq A_k\!\sum_{i=0}^k\sum_{n=1}^\infty
\frac{1}{n^2}
\underbrace{n^{4k+2}\|\gamma_n^{<i>}\|_{n+2k,\infty}}_{\leq B_k}
\leq A_kB_k\!\sum_{i=0}^k\sum_{n=1}^\infty\frac{1}{n^2}
< \infty .\vspace{-1.9mm}
\]
Thus $\sum_{n=1}^\infty\eta_n$
is absolutely convergent
and hence convergent in $BC^\infty(\R,\wb{E})$.
Pointwise calculation of the limit shows that
$\sum_{n=1}^\infty\eta_n=\gamma$
from above.\linebreak
Since $t_\infty\not\in\Supp(\eta_n)$
for all $n\in \N$, we can argue now
as at the end of the proof of Theorem~\ref{thmgencu}
to see that
$\gamma\in BC^\infty(\R,E)$.
By construction, $\gamma$ has also all
other required properties.\,\vspace{-2mm}\Punkt
\appendix
\section{Proofs of the lemmas in Section~\ref{secperlim}}
In this appendix,
proofs are provided for the lemmas from Section~\ref{secperlim}.\\[2.5mm]
{\bf Proof of Lemma~\ref{completeBC}.}
(a) For each $0$-neighbourhood~$U$ in~$E$,
we set $\lfloor X,U\rfloor:=\{\gamma\in BC(X,E)\colon
\gamma(X)\sub U\}$.
If $V\sub E$ is a $0$-neighbourhood such that
$V= -V$ and $V+V\sub U$,
then
$\lfloor X,V\rfloor+\lfloor X,V\rfloor \sub \lfloor X,U\rfloor$
and $-\lfloor X,V\rfloor\sub \lfloor X,U\rfloor$,
entailing that there is a unique group
topology on $BC(X,E)$ for which the sets
$\lfloor X,U\rfloor$ form a basis
of $0$-neighbourhoods. As $\bigcap_U U=\{0\}$,
also the sets $\lfloor X,U\rfloor$
have intersection $\{0\}$ and thus $BC(X,E)$
is Hausdorff.
To see that the given group topology is a vector topology,
it only remains to check conditions
$(\mbox{EVT}_{\text{I}}')$--$(\mbox{EVT}_{\text{III}}')$
of \cite[Ch.\,I, \S1, no.\,1]{BTV}.
First, given $\gamma_0\in BC(X,E)$,
we show that the map $\K\to BC(X,E)$, $t\mto t\gamma_0$
is continuous at~$0$.
Since $\gamma_0(X)$ is bounded, for each $0$-neighbourhood
$U\sub E$ there is a $0$-neighbourhood
$V\sub \K$ such that $V\gamma_0(X)\sub U$.
Then $V\gamma_0\sub \lfloor X,U\rfloor$,
entailing the assertion.
Next, given $t_0\in \K$, let us check that the map
$BC(X,E)\to BC(X,E)$, $\gamma\mto t_0\gamma$
is continuous at~$0$.
In fact, given~$U$ as before,
there is a $0$-neighbourhood $V\sub E$ such that
$t_0V\sub U$. Then $t_0\lfloor X,V\rfloor\sub \lfloor X,U\rfloor$.
Also, scalar multiplication
$\K\times BC(X,E)\to BC(X,E)$,
$(t,\gamma)\mto t\gamma$
is continuous at~$(0,0)$.
In fact, given~$U$,
there are $0$-neighbourhoods
$V\sub \K$ and $W\sub E$ such that
$VW\sub U$. Then $V\lfloor X,W\rfloor \sub
\lfloor X,U\rfloor$.\\[3mm]
(b) If $(\gamma_\alpha)_\alpha$ is a
Cauchy net in $BC(X,E)$,
then $(\gamma_\alpha(x))_\alpha$
is a Cauchy net in~$E$
for each $x\in X$,
the point evaluation $BC(X,E)\to E$, $\gamma\mto\gamma(x)$
being continuous linear.
Since $E$ is complete, $\gamma_\alpha(x)\to\gamma(x)$
for some $\gamma(x)\in E$.
Given a $0$-neighbourhood $U\sub E$,
let $V\sub E$ be a $0$-neighbourhood
such that $V+V+V\sub U$,
and $W\sub V$ be a closed, symmetric
$0$-neighbourhood
such that $SW\sub V$ for some $0$-neighbourhood
$S\sub \K$.
There exists $\alpha_0$ such that
$\gamma_\alpha-\gamma_\beta\in \lfloor X,W\rfloor$
for all $\alpha,\beta\geq \alpha_0$.
Then $\gamma_\alpha(x)-\gamma_\beta(x)\in W$ for each $x\in X$.
Since~$W$ is closed, passage to the limit yields
%
\begin{equation}\label{give3}
\gamma_\alpha(x)-\gamma(x)\,\in \, W \,\sub \, V \quad\mbox{for each $x\in X$
and $\alpha\geq \alpha_0$.}
\end{equation}
Each $x_0\in X$ has a neighbourhood~$Q$
such that $\gamma_{\alpha_0}(x)
-\gamma_{\alpha_0}(x_0)\in V$ for all $x\in Q$
and hence $\gamma(x)-\gamma(x_0)=
(\gamma(x)-\gamma_{\alpha_0}(x))
+(\gamma_{\alpha_0}(x)-\gamma_{\alpha_0}(x_0))$
$+(\gamma_{\alpha_0}(x_0)-\gamma(x_0))\in U$.
Thus $\gamma$ is continuous at~$x_0$ and
hence continuous.
To see that~$\gamma(X)$ is bounded,
let $T\sub S$ be a $0$-neighbourhood
such that $T\gamma_{\alpha_0}(X)\sub V$.
Then
$t\gamma(x)=t\gamma_{\alpha_0}(x)+t(\gamma(x)-\gamma_{\alpha_0}(x))
\in V+SW\sub U$
for each $x\in X$ and $t\in T$,
whence $T\gamma(X)\sub U$.
Thus $\gamma\in BC(X,E)$.
Since $\gamma_\alpha-\gamma\in \lfloor X,U\rfloor$
for all $\alpha\geq \alpha_0$,
we see that $\gamma_\alpha\to\gamma$.\\[3mm]
(c) It is obvious that $\theta$ is linear,
and it is a topological embedding
by definition of the topology on $BC^k(U,E)$.
To see that the image is closed,
let $(\gamma_\alpha)_\alpha$ be a net in
$BC^k(U,E)$ such that $\theta(\gamma_\alpha)\to\eta$
for some $\eta=(\eta_j)\in \prod_j BC(U^{j+1},E)$.
We claim that $\gamma:=\eta_0\in BC^k(U,E)$
and $\theta(\gamma)=\eta$.
This will be the case if
%
\begin{equation}\label{notn}
\eta_{j+1}(x_0,\ldots,x_{j+1})\,=\,
\frac{\eta_j(x_0,x_1,\ldots,x_j)-\eta_j(x_{j+1},x_1,\ldots,
x_j)}{x_0-x_{j+1}}
\end{equation}
for each $j\in \N_0$ with $j<k$
and each $(x_0,\ldots,x_{j+1})\in U^{j+2}$
with $x_0\not=x_{j+1}$.
To prove (\ref{notn}),
we use that the $j$th component $\theta_j(\gamma_\alpha)$
converges to~$\eta_j$,
and the continuity of the
point evaluation
$\varepsilon_1
\colon
BC(U^{j+1},E)\to E$,
$\zeta\mapsto \zeta(x_0,\ldots,x_j)$,
the point evaluation
$\varepsilon_2\colon
BC(U^{j+1},E)\to E$ at $(x_{j+1},x_1,\ldots, x_j)$
and the point evaluation
$\varepsilon_3
\colon
BC(U^{j+2},E)\to E$ at $(x_0,\ldots, x_{j+1})$.
Since $\varepsilon_3(\theta_{j+1}(\gamma_\alpha))
=\gamma_\alpha^{<j+1>}(x_0,\ldots, x_{j+1})
=\frac{\varepsilon_1(\theta_j(\gamma_\alpha))
-\varepsilon_2(\theta_j(\gamma_\alpha))}{x_0-x_{j+1}}$,
passing to the limit we obtain~(\ref{notn}).\\[3mm]
(d) It is clear from (c) that $BC^k(U,E)$ is
a topological $\K$-vector space.
If~$E$ is complete, then also $BC^k(U,E)$
is complete as it is isomorphic to a closed
vector subspace of a complete
topological vector space by (b) and (c).\vspace{2.7mm}\Punkt

\noindent
{\bf Proof of Lemma~\ref{spla1}.}
It is clear that the map $\gamma\mto\|\gamma^{<j>}\|_{q,\infty}$
is positively homogeneous.
Given $r>0$, pick $s\in \; ]0,r[$.
Then the set $V:=B_s^q(0)$
is a $0$-neighbourhood in~$E$ and
hence $\lfloor U^{j+1},V\rfloor$
is $0$-neighbourhood
in $BC(U^{j+1},E)$,
entailing that $W:=\{\gamma\in BC^k(U,E)\colon
\gamma^{<j>}\in \lfloor U^{j+1},V\rfloor\}$
is a $0$-neighbourhood in $BC^k(U,E)$.
As $W\sub
\{\gamma\in BC^k(U,E)\colon \|\gamma^{<j>}\|_{q,\infty}\leq s\}$,
we see that also
$\{\gamma\in BC^k(U,E)\colon \|\gamma^{<j>}\|_{q,\infty}<r\}$
is a $0$-neighbourhood.
Hence the mappings in contention are
gauges.\\[2.7mm]
To see that a fundamental
system of gauges is obtained,
we use that each
$0$-neighbourhood in $BC^k(U,E)$
contains a finite intersection
of sets of the form
\[
W\; :=\; \{\gamma\in BC^k(U,E)\colon
\gamma^{<j>}\in \lfloor U^{j+1},V\rfloor\}\,,
\]
where $j\in \N_0$ with $j\leq k$
and $V\sub E$ is a $0$-neighbourhood.
There are $r_1,\ldots, r_n>0$ and $q_1,\ldots, q_n\in \Gamma$
such that $\bigcap_{i=1}^n B_{r_i}^{q_i}(0)\sub V$.
Consider the gauges
$p_i\colon \gamma\mto \|\gamma^{<j>}\|_{q_i,\infty}$
for $i\in \{1,\ldots, n\}$.
Then $\bigcap_{i=1}^n B^{p_i}_{r_i}(0)
\sub W$.\vspace{2.7mm}\Punkt

\noindent
{\bf Proof of Lemma~\ref{transhomo}.}
(a)
A trivial induction on $j\in \{0,1,\ldots, k\}$
gives $\eta^{<j>}(x_0,\ldots,x_j)=
\gamma^{<j>}(x_0+t_0,\ldots, x_j+t_0)$
for all $(x_0,\ldots, x_j)\in (U-t_0)^{j+1}$.
Now $\|\eta^{<k>}\|_{q,\infty}=\|\gamma^{<k>}\|_{q,\infty}$
is an immediate consequence.\\[2.5mm]
(b) A trivial induction on $j\in \{0,1,\ldots, k\}$
shows that $\eta^{<j>}(x_0,\ldots,x_j)=
a^j\gamma^{<j>}(ax_0,\ldots, ax_j)$
for all $(x_0,\ldots, x_j)\in (a^{-1}U)^{j+1}$.
Now $\|\eta^{<k>}\|_{q,\infty}=
|a|^k \|\gamma^{<k>}\|_{q,\infty}$
is an immediate consequence.\\[2.5mm]
(c) A trivial induction on $j\in \{0,1,\ldots, k\}$
shows that $\gamma|_V\in BC^j(V,E)$
and $(\gamma|_V)^{<j>}=\gamma^{<j>}|_{V^{j+1}}$.
Now $\|(\gamma|_V)^{<k>}\|_{q,\infty}\leq
\|\gamma^{<k>}\|_{q,\infty}$
is immediate.\,\vspace{2.7mm}\Punkt

\noindent
{\bf Proof of Lemma~\ref{charctbd}.}
If $\gamma$ is $C^k$ in our sense,
then it is easy to show by induction
on~$j$ that $\gamma$ is $C^j$ in the usual
sense for each $j\in \{0,\ldots, k\}$, with
\begin{equation}\label{trivboo}
\gamma^{(j)}(x)\; =\; j! \gamma^{<j>}(x,\ldots,x)\quad
\mbox{for all $x\in I$}
\end{equation}
(cf.\ \cite[Proposition~6.2]{BGN} and \cite[\S29]{Sch}).
Furthermore,
it is clear from the preceding formula
that $\gamma^{(j)}(I)$
is bounded if so is $\gamma^{<j>}(I^{j+1})$.\\[2.7mm]
Conversely, assume that $\gamma$
is continuous and assume that the derivatives
$\gamma',\gamma'',\ldots,\gamma^{(k)}$
exist and are continuous.
Then we have, for all $x_0\not=x_1$ in~$I$,
\[
\frac{\gamma(x_1)-\gamma(x_0)}{x_1-x_0}
\,=\,\int_{t_1=0}^1 \gamma'(x_0+t_1(x_1-x_0))\, dt_1
\]
by the Fundamental Theorem
of Calculus for curves in locally convex
spaces (see, e.g., \cite[Chapter~1]{GaN}).
Hence
\[
\gamma^{<1>}\colon I^2\to E\,,\quad
(x_0,x_1)\mto
\int_{t_1=0}^1 \gamma'((1-t_1)x_0+t_1x_1)\, dt_1
\]
is an extension to the difference
quotient map, which is
continuous by
the theorem on parameter-dependence of weak integrals
(see \cite[Chapter~1]{GaN}.
Thus~$\gamma$ is~$C^1$ in the sense of
Definition~\ref{def1}.
Iterating the argument, we find
that $\gamma$ is $C^j$ for each $j\in \{1,\ldots, k\}$
and $\gamma^{<j>}(x_0,x_1,\ldots, x_j)$
is given by
%
%
\begin{eqnarray}
\lefteqn{\int_{t_1=0}^1 \cdots \int_{t_j=0}^1
(1-t_{j-1})\cdots (1-t_1)\cdot}\qquad  \label{enabothest}\\
& & \quad\;\; \gamma^{(j)}(t_1x_1 + (1-t_1)t_2x_2 + \cdots+
(1-t_1)(1-t_2)\cdots(1-t_{j-1})t_jx_j\notag \\[1mm]
& & \qquad\qquad  +\,(1-t_1)\cdots(1-t_j)x_0)\;dt_j\ldots dt_1\,.\notag
\end{eqnarray}
If $\gamma^{(j)}$ is bounded, then so is $\gamma^{<j>}$
by the preceding formula, with
%
\begin{equation}\label{oppoesti}
\|\gamma^{<j>}\|_{q,\infty}\;\leq\;
\|\gamma^{(j)}\|_{q,\infty}
\end{equation}
for each continuous seminorm~$q$ on~$E$.
This completes the proof.\,\vspace{5mm}\Punkt

\noindent
{\bf Remark.}
By (\ref{trivboo})
and (\ref{oppoesti}),
the seminorms $\gamma\mto
\|\gamma^{<j>}\|_{q,\infty}$
(for $j\leq k$ and $q$ in
the set of continuous seminorms
on~$E$) define the same vector topology
on $BC^k(U,E)$ as the seminorms
$\gamma\mto\|\gamma^{(j)}\|_{q,\infty}$
ordinarily used on
this space.\\[4mm]
\emph{Acknowledgement.}
The author thanks S.\,V. Ludkovsky for
remarks in the final stages of the preparation
of this article.\vspace{-2mm}
{\small
{\bf Helge Gl\"{o}ckner},
TU Darmstadt, FB Mathematik AG~5,
Schlossgartenstr.\,7,\\
64289 Darmstadt, Germany; FAX: +49-6151-166030\\[.4mm]
\,E-mail: {\tt gloeckner@mathematik.tu-darmstadt.de}}

\begin{thebibliography}{99}\itemsep+1.6pt
%
%
\bibitem{AEK}
Adasch, N., B. Ernst and D. Keim,
``Topological Vector Spaces.
The Theory without Convexity Conditions,''
Springer, 1978.
%
%
\bibitem{Ber} Bertram, W.,
\emph{Differential geometry, Lie groups
and symmetric spaces over general base fields and rings},
to appear in Memoirs of the AMS
(cf.\ arXiv:math/0502168).
%
%
\bibitem{BGN}
Bertram, W., H. Gl\"{o}ckner and K.-H. Neeb,
\emph{Differential calculus over general base fields and rings},
Expo.\ Math.\ \textbf{22} (2004),
213--282.
%
%
\bibitem{Bom}
Boman, J.
\emph{Differentiability of a function and of
its compositions with functions of one variable},
Math Scand.\ \textbf{20} (1967), 249--268.
%
%
\bibitem{BTV}
Bourbaki, N., ``Topological Vector Spaces''
Chapters~1--5, Springer-Verlag, 1987.
%
%
\bibitem{DSm}
De Smedt, S.
\emph{$p$-adic continuously differentiable functions of
several variables}, Collect.\ Math.\ \textbf{45}
(1994), 137--152.
%
%
\bibitem{FaK}
Fr\"{o}licher, A. and A. Kriegl,
``Linear Spaces and Differentiation Theory,''
John Wiley, 1988.
%
%
\bibitem{NOA}
Gl\"{o}ckner, H.,
\emph{Smooth Lie groups over local fields of positive characteristic
need not be analytic}, J. Algebra
\textbf{285} (2005), 356--371.
%
%
\bibitem{ANA}
Gl\"{o}ckner, H.,
\emph{Every smooth $p$-adic Lie group admits a
compatible analytic structure},
Forum Math.\ {\bf 18} (2006), 45--84.
%
%
\bibitem{SUR}
Gl\"{o}ckner, H.,
\emph{Aspects of $p$-adic non-linear functional analysis},
pp.\ 237--253 in:
A.\,Yu.\ Khrennikov, Z. Raki\'{c}
and I.\,V. Volovich (Eds.),
$p$-Adic
Mathematical Physics.
2nd International Conference (Belgrade, 2005),
AIP Conference Proceedings \textbf{826},
Amer.\ Inst.\ Physics, New York, 2006
(cf.\ arXiv:math/0602081).
%
%
\bibitem{ZOO}
Gl\"{o}ckner, H.,
\emph{Lie groups over non-discrete topological fields},
preprint, arXiv:math/0408008\,.
%
%
\bibitem{IM2}
Gl\"{o}ckner, H.,
\emph{Finite order differentiability properties, fixed points
and implicit functions over valued fields},
preprint, arXiv:math/0511218\,.
%
%
\bibitem{COM}
Gl\"{o}ckner, H., \emph{Comparison of some
notions of $C^k$-maps in multi-variable
non-archimedian analysis},
preprint, arXiv:math/0609041\,.
%
%
\bibitem{GaN}
Gl\"{o}ckner, H. and K.-H. Neeb,
``Infinite-Dimensional Lie Groups,'' Vol.\,I;
book in preparation.
%
%
\bibitem{Jar}
Jarchow, H., ``Locally Convex Spaces,''
B.\,G. Teubner, Stuttgart, 1981.
%
%
\bibitem{Kal}
Kalton, N.\,J., N.\,T. Peck and J.\,W. Roberts,
``An F-space Sampler,''
Cambridge University Press, Cambridge,
1984.
%
%
\bibitem{KaM}
Kriegl, A. and P.\,W. Michor,
``The Convenient Setting of Global
Analysis,'' Amer.\ Math.\ Soc., Providence, 1997.
%
%
\bibitem{Ld1}
Ludkovsky, S.\,V.,
\emph{Irreducible unitary representations
of non-archimedian groups of diffeomorphisms},
Southeast Asian Bull.\ Math.\
\textbf{22} (1998), 419--436.
%
%
\bibitem{Ld2}
Ludkovsky, S.\,V.,
\emph{Quasi-invariant measures on
non-Archimedian groups and
semigroups of loops and paths,
their representations} I, II
Ann.\ Math.\ Blaise Pascal \textbf{7} (2000),
19--53 and 55--80.
%
%
\bibitem{Lud}
Ludkovsky, S.\,V.,
\emph{Smoothness of functions
global and along curves over
ultra-metric fields},
arXiv:math/0608725\,.
%
%
\bibitem{Sch}
Schikhof, W.\,H.,
``Ultrametric Calculus,''
Cambridge University Press,
1984.\vspace{1.7mm}
%
\end{thebibliography}
\end{document}